\documentclass[a4paper,10pt]{article}
\usepackage[psamsfonts]{amssymb}
\usepackage[psamsfonts]{eucal}

\usepackage{amsmath}
\usepackage{theorem}
\newtheorem{thm}{Theorem}

\newtheorem{cor}[thm]{Corollary}

\theorembodyfont{\rmfamily}

\theorembodyfont{\rmfamily}

\begin{document}

\title{{GEOMETRY OF NEUTRAL METRICS IN DIMENSION FOUR}
}
\author{{Johann Davidov, Gueo Grantcharov and Oleg Mushkarov \thanks{Partially supported by CNRS-BAS joint research project {\it
Invariant metrics and complex geometry, 2008-2009.}}}}

\date{}
\maketitle

 \rm
\begin{abstract}
 The purpose of this article is to review some recent results on the geometry
of neutral signature metrics in dimension four and their twistor spaces. The
following topics are considered: Neutral K\"ahler and hyperk\"ahler surfaces,
Walker metrics, Neutral anti-self-dual 4-manifolds and projective structures,
Twistor spaces of neutral metrics.



\end{abstract}

\vspace{0.5cm}

\section{\bf Introduction}

 Riemannian and Lorentzian metrics in dimension four have been studied extensively for a
long time in connection with various important problems in geometry and
physics. The efforts were concentrated mainly on understanding the
corresponding Einstein equations, which led to deep results in both fields. The
metrics of neutral signature $(+,+,-,-)$ appear also in many geometric and
physics problems but they have received less attention until recently. An
impulse for the development of the geometry in this signature was the work of
Ooguri and Vafa \cite{OV} who showed that this setting naturally arises in
$N=2$ string theory. Moreover, it was realized that powerful techniques from
the Seiberg-Witten theory and the theory of integrable systems can be
successfully used to study K\"ahler-Einstein and self-dual metrics as well as
the self-dual Yang-Mills equations in neutral signature (see \cite{P},\cite{DW}
and \cite{Uhl} and references therein).

At the linear level, there are close analogues between the neutral signature
and Riemannian metrics. For instance, in both cases the Hodge star operator is
an involution of the bundle of two-forms, which leads to a splitting of the
curvature operator into four irreducible components. Thus the most significant
classes of Riemannian metrics like Einstein and (anti) self-dual have natural
counterparts in the neutral case. However, there are profound differences
between the two geometries, both locally and globally.  For example, in the
neutral case the gauged-fixed self-duality equations are ultrahyperbolic,
whereas in the Riemannian case they are elliptic, which suggests that the
neutral case is less rigid than the Riemannian one. Indeed, any Riemannian
self-dual conformal structure must be real-analytic but this is not true for
neutral structures.

 In the present survey, we review some recent results on the geometry of neutral signature
 metrics on 4-manifolds with emphasis on the differences with the Riemannian
 case.

  In Sections 3 and 4, we present the results of Petean \cite{P} and Kamada \cite{Ka}
 about the  existence of neutral K\"ahler-Einstein and neutral hyperk\"ahler metrics
 on compact complex surfaces.

 Section 5 is devoted to the so-called Walker metrics which have interesting curvature
 properties and provide examples showing that some local integrability results related to
 the Goldberg conjecture \cite{Go}
 are not true for neutral metrics.

 In Section 6, we describe the Dunajski-West \cite{DW1} correspondence between
the two dimensional projective structures and the anti-self-dual neutral 4-manifolds with a non-trivial null
conformal Killing vector field.

Finally, in Section 7, we review some geometric and analytic results on the
hyperbolic twistor spaces \cite{BDM} and discuss their relation to the
LeBrun-Mason twistor spaces of neutral self-dual manifolds \cite{LB-M}.

This survey on the geometry of neutral metrics in dimension four is by no means comprehensive and it is in a
sense complementary to \cite{DW}, where the connection of neutral anti-self-dual conformal structures with
integrable systems and twistor theory is underlined.

It is a pleasure to thank the Program Committee of the Thirty Seventh Spring
Conference of the Union of Bulgarian Mathematicians for the kind invitation to
deliver a lecture on geometry of neutral metrics.

\section{\bf Preliminaries}

 A pseudo-Riemannian metric on a smooth 4-manifold $M$ is called $neutral$ if it has signature
$(+,+,-,-)$. In contrast to the Riemannian case the existence of such metrics
on compact manifolds imposes topological restrictions since it is equivalent to
the existence of a field of 2-planes. We refer to
\cite{Matsushita1991,LawMatsushita,DW} for more information.

   Let $M$ be an oriented 4-manifold with a
neutral metric $g$. Then $g$ induces an inner product on the bundle $\Lambda^2$
of bivectors by
$$
<X_1\wedge X_2,X_3\wedge
X_4>=\frac{1}{2}[g(X_1,X_3)g(X_2,X_4)-g(X_1,X_4)g(X_2,X_3)],
$$
$X_1,...,X_4\in TM$. Let ${\bf e}_1,\ldots,{\bf e}_4$ be a local oriented
orthonormal frame of $TM$ with $||{\bf e}_1||^2=||{\bf e}_2||^2=1$, $||{\bf
e}_3||^2=||{\bf e}_4||^2=-1$. As in the Riemannian case, the Hodge star
operator $\ast:\Lambda^2\to \Lambda^2$ is an involution given by
$$*({\bf e}_1\wedge{\bf e}_2)={\bf e}_3\wedge{\bf
e}_4,\quad *({\bf e}_1\wedge{\bf e}_3)={\bf e}_2\wedge{\bf e}_4,\quad *({\bf
e}_1\wedge{\bf e}_4)=-{\bf e}_2\wedge{\bf e}_3.$$
 Denote by $\Lambda_{\pm}$ the
subbundles of $\Lambda^2$ determined by the eigenvalues $\pm 1$ of the Hodge
star operator. Set
\begin{equation}\label{eq s-s-basis}
\begin{array}{lll}
s_1={\bf e}_1\wedge{\bf e}_2-{\bf e}_3\wedge{\bf e}_4,\quad\quad \bar s_1={\bf
e}_1\wedge{\bf e}_2+{\bf e}_3\wedge{\bf e}_4, \\
s_2={\bf e}_1\wedge{\bf e}_3-{\bf e}_2\wedge{\bf e}_4,\quad\quad \bar s_2={\bf
e}_1\wedge{\bf e}_3+{\bf e}_2\wedge{\bf e}_4, \\
s_3={\bf e}_1\wedge{\bf e}_4+{\bf e}_2\wedge{\bf e}_3,\quad\quad \bar s_3={\bf
e}_1\wedge{\bf e}_4-{\bf e}_2\wedge{\bf e}_3.
\end{array}
\end{equation}
Then $\{s_1,s_2,s_3\}$ and $\{\bar s_1,\bar s_2,\bar s_3\}$ are local oriented
orthonormal frames of $\Lambda_{-}$ and $\Lambda_{+}$ respectively with
$||s_1||^2=||\bar s_1||^2=1$, $||s_2||^2=||\bar s_2||^2=||s_3||^2=||\bar
s_3||^2=-1$.

   Let ${\cal R}:\Lambda^2 \longrightarrow\Lambda^2$ be the curvature operator of $(M,g)$.
It is related to the curvature tensor $R$ by
$$g({\cal R}(X\wedge Y),Z\wedge T)=g(R(X,Y)Z,T);\quad
X,Y,Z,T\in TM.$$ In this paper, we adopt the following definition of the curvature tensor
$R(X,Y)=\nabla_{[X,Y]}-[\nabla_{X},\nabla_{Y}]$, where $\nabla$ is the Levi-Civita connection of $g$. The
curvature operator ${\cal R}$ admits an $SO(2,2)$-irreducible decomposition
$${\cal R}={\tau\over 6}I+{\cal B}+{\cal W}_{+}+{\cal W}_{-}$$
similar to that in the 4-dimensional Riemannian case \cite{bes}. Here $\tau$ is
the scalar curvature, $\cal B$ represents the traceless Ricci tensor, ${\cal
W}={\cal W}_{+}+{\cal W}_{-}$ corresponds to the Weyl conformal tensor, and
${\cal W}_{\pm}={\cal W}|\Lambda_{\pm} = \displaystyle{\frac{1}{2}({\cal
W}\pm\ast{\cal W})}$. The metric $g$ is Einstein exactly when ${\cal B}=0$ and
is conformally flat when ${\cal W}=0$. It is said to be self-dual, resp.
anti-self-dual,  if ${\cal W}_{-}=0$, resp. ${\cal W}_{+}=0$.

 A neutral almost Hermitian structure on $M$ consists of
an almost complex structure $J$ and a neutral metric $g$ satisfying the compatibility condition
$g(JX,JY)=g(X,Y)$. If the almost complex structure $J$ is integrable (i.e. comes from a complex structure on
$M$), the structure $(g,J)$ is called Hermitian. As in the positive definite case, any neutral almost Hermitian
structure $(g,J)$ determines a non-degenerate $2$-form $\Omega(X,Y)=g(JX,Y)$ but now $\Omega$ is compatible with
the opposite orientation of $M$. If $\Omega$ is closed (i.e., it is a symplectic form) the structure is said to
be neutral almost K\"{a}hler and it is called neutral K\"{a}hler if, in addition, the almost complex structure
$J$ is integrable.

  To describe the neutral almost Hermitian structures on $(M,g)$ in terms of
 bivectors,  we identify $\Lambda^2$
with the bundle of skew-symmetric endomorphisms of $TM$ by assigning to each
$\sigma\in\Lambda^2$ the endomorphism $J_{\sigma}$ on $T_pM$, $p=\pi(\sigma)$,
defined by
\begin{equation}\label{3.1}
g(J_{\sigma}X,Y)=2g(\sigma,X\wedge Y); X,Y\in T_pM.
\end{equation}
Then the local sections $s_1,s_2,s_3$ of $\Lambda_-$, defined by (\ref{eq s-s-basis}) determine local
endomorphisms $J_1,J_2,J_3$ of $TM$ that satisfy the relations
\begin{equation}\label{3.2}
J_1^2=-J_2^2=-J_3^2=-Id, \hspace{.1in} J_1J_2=-J_2J_1=J_3
\end{equation}
of the imaginary units of the paraquaternionic algebra (split quaternions).
Hence the set of all
$$J = y_1J_1+y_2J_2+y_3J_3, \hspace{.1in} y_1^2-y_2^2-y_3^2= 1$$ describes the
local almost complex structures compatible with $g$ and the orientation of $M$.

 If there are three global integrable structures $J_1, J_2, J_3$ on $M$
satisfying the relations (\ref{3.2}), then the structure $(g,J_1, J_2, J_3)$ is called neutral hyperhermitian.
When additionally the 2-forms $\Omega_i(X,Y) = g(J_iX,Y)$ are closed, the hyperhermitian structure is called
neutral hyperk\"ahler. It is well known \cite{Kamada} that the neutral hyperhermitian metrics are self-dual,
whereas the neutral hyperk\"ahler metrics are self-dual and Ricci-flat.



\vspace{0.5cm}

\section{Neutral K\"ahler surfaces}

It is well known that a compact complex surface admits a K\"ahler metric if and
only if its first Betti number is even (see e.g. \cite{BPV}). The analogous
problem for neutral K\"ahler metrics is not yet solved completely. Using the
Seiberg-Witten theory Petean \cite{P} proved the following:

\begin{thm}
Let $(M,g,J)$ be a compact neutral K\"ahler surface and $k(M,J)$ its Kodaira
dimension.

(i) If $k(M,J)=-\infty$, then $(M,J)$ is either a ruled surface, or a surface
of class $VII_0$ with no global spherical shell and with positive even second
Betti number.

(ii) If $k(M,J)=0$, then $(M,J)$ is either a hyperelliptic surface, a primary
Kodaira surface or a complex torus.

(iii) If $k(M,J)=1$, then $(M,J)$ is a minimal properly elliptic surface with
zero signature.

(iv) If $k(M,J)=2$, then (M,J) is a minimal surface of general type with non-negative even signature.

\end{thm}

 As it is noted in \cite{Kamada}, all surfaces listed in (ii) and (iii) as well as the
ruled surfaces in (i) do admit neutral K\"ahler metrics but the existence of such metrics in the remaining cases
is still open. Moreover, there is a conjecture by Nakamura \cite{Nak}, which states that there are no surfaces
satisfying the conditions in the second case of (i).

Now we turn to the existence of self-dual neutral K\"ahler metrics noting that
in this case the self-duality is equivalent to the vanishing of the scalar
curvature. The following result of  Kamada \cite{Kamada} was obtained by means
of Theorem 1.

\medskip

\begin{thm}
Let $(M,g,J)$ be a compact self-dual neutral K\"ahler surface and let $\tau(M)$
and $c_1(M,J)$ be its signature and first Chern class, respectively. Then
$(M,J)$ is one of the following:

(i) A surface of class $VII_0$ with no global spherical shell and with positive
even second Betti number,  if $c_1(M,J)^2<0$.

(ii) A minimal surface of general type with even signature,  if $c_1(M,J)^2>0$
and $\tau(M)>0$.

(iii) A Hirzebruch surface or a minimal surface of general type uniformized by
the bidisc,  if $c_1(M,J)^2>0$ and $\tau(M)=0$.

(iv) A hyperelliptic surface, a primary Kodaira surface, a complex torus or a
minimal properly elliptic surface with zero signature, if $c_1(M,J)^2=0$.

\end{thm}

 A typical example of compact self-dual neutral K\"ahler surface is the
product ${\Bbb C}{\Bbb P}^1\times{\Bbb C}{\Bbb P}^1$ with the neutral product metric $g=h\oplus -h$, where $h$
is the round metric on ${\Bbb C}{\Bbb P}^1=S^2$. Note that this metric is not only self-dual but also
conformally flat and by a result of Kuiper \cite{kuiper} any compact simply-connected conformally flat neutral
4-manifold is equivalent to $(S^2\times S^2, g)$. On the other hand, Tod \cite{tod} constructed deformations of
$g$ consisting of $S^1$-invariant self-dual neutral K\"ahler metrics and recently Kamada \cite{kamada2} proved
that such metrics with a Hamiltonian $S^1$-symmetry can exist only on $S^2\times S^2$.

We should note also \cite{Kamada} that it is not known whether the surfaces in Theorem 2\,(ii) and the properly
elliptic surfaces with zero signature admit self-dual neutral K\"ahler metrics.

 The following theorem of Petean \cite{P} gives a complete classification of
compact complex surfaces admitting neutral  Ricci-flat K\"ahler metrics as well
as an almost complete classification of surfaces that admit neutral
K\"ahler-Einstein metrics with non-zero scalar curvature.

\begin{thm}\label{petean}

Let $(M,g,J)$ be a compact neutral K\"ahler-Einstein surface. Then  $(M,J)$ is
one of the following:

(i) A complex torus.

(ii) A hyperelliptic surface.

(iii) A primary Kodaira surface.

(iv) A minimal ruled surface over a curve of genus ${\bf g}\geq 2$.

(v) A minimal surface of class $VII_0$ with no global spherical shell, and with
even and positive second Betti number.

If $(M,J)$ admits a neutral Ricci-flat K\"ahler metric, then it is as in (i),
(ii) or (iii).

\end{thm}

Petean \cite{P} constructed also examples of neutral K\"ahler - Einstein metrics on the surfaces of type $(i),
(ii), (iii)$ and on most surfaces of type $(iv)$. We note that the Petean's metrics on a complex torus are
actually hyperk\"ahler and depend on an arbitrary positive smooth function on an elliptic curve. This shows that
the moduli spaces of neutral Ricci-flat K\"ahler metrics can be highly non-trivial and different from those in
the positive definite case.

\section{Neutral hyperk\"ahler surfaces}

As is well known (c.f. \cite{bes}), any compact hyperk\"ahler surface is either a complex torus with a flat
metric or a $K3$-surface with Calabi-Yau metric. By contrast, H.Kamada \cite{Ka} proved that neutral
hyperk\"ahler structures can exist either on a complex torus or on a primary Kodaira surface and obtained a
description of all such structures on both types of surfaces. Here we shall provide his result for primary
Kodaira surfaces.

 Consider the affine transformations
  $\rho_i(z_1,z_2) = (z_1+a_i,z_2+\overline{a}_iz_1+b_i)$ of ${\Bbb C}^2$,
 where $a_i$,$b _i$, $i=1,2,3,4$, are  complex
numbers such that $a_1=a_2=0, Im(a_3{\overline a}_4) = b_1$. Then $\rho_i$
generate a group $G$ of affine transformations acting freely and properly
discontinuously on ${\Bbb C}^2$. The quotient space ${\Bbb C}^2/G$ is called a
primary Kodaira surface.

 Note also \cite{Hit1990,Ka} that any neutral hyperk\"ahler structure is
 determined by three symplectic forms
 $(\Omega_1, \Omega_2 , \Omega_3)$ satisfying the relations
 $$-\Omega_1^2=\Omega_2^2=\Omega_3^2, \quad \Omega_l\wedge\Omega_m = 0,\> l\neq m.$$
 In these terms the result of Kamada \cite{Ka} is the following:

 \begin{thm}
For any neutral hyperk\"ahler structure on a primary Kodaira surface $M$ there
are complex coordinates $(z_1, z_2)$ of ${\Bbb C}^2$ such that the structure is
given by the following symplectic forms:
$$\Omega_1 = Im(dz_1\wedge d\overline{z}_2)+iRe(z_1)dz_1\wedge d\overline{z}_1 +
(i/2)\partial\overline{\partial}\phi, \\
$$
$$
\Omega_2 = Re(e^{i\theta}dz_1\wedge d\overline{z}_2),\hspace{.1in} \Omega_3 =
Im(e^{i\theta}dz_1\wedge d\overline{z}_2),
$$

where $\theta$ is a real constant and $\phi$ is a smooth function on $M$ such
that:

\begin{equation}\label{kod-hk}
4i(Im(dz_1\wedge d\overline{z}_2) +iRe(z_1)(dz_1\wedge
d\overline{z}_1))\wedge\partial\overline{\partial}\phi =
\partial\overline{\partial}\phi\wedge\partial\overline{\partial}\phi
\end{equation}

 \end{thm}

As proved by Kamada \cite{Ka}, the neutral hyperk\"ahler metric on $M$ defined by $(\Omega_1, \Omega_2 ,
\Omega_3)$ is flat if and only if the function $\phi$ is constant. On the other side, any primary Kodaira
surface is a toric bundle over an elliptic curve and the pull-back of any smooth function on the base curve
gives a solution to (\ref{kod-hk}). It follows that the moduli space of neutral hyperk\"ahler structures on a
primary Kodaira surface is infinite dimensional, which is in sharp contrast with the positive definite case. We
refer to \cite{Kamada} for analogous description of neutral hyperk\"ahler structures on complex tori.

\section{\bf Walker metrics}

A basic problem in almost Hermitian geometry is to relate  properties of an
almost Hermitian structure $(g,J)$ to the curvature of the metric $g$. For
example, the well-known Goldberg conjecture \cite{Go} claims that a compact
almost K\"{a}hler manifold is K\"ahler  provided the metric $g$ is Einstein.
This conjecture was proved by Sekigawa \cite{Sek} in the case of non-negative
scalar curvature but it is still far from being solved in the negative case. We
refer to the survey \cite{Apo-Dra2003} for an update on the integrability of
almost K\"{a}hler structures. Another integrability result related to the
curvature properties of a manifold is the Riemannian version of the well-known
Goldberg-Sacks theorem in General Relativity. It says that an oriented Einstein
$4$-manifold admits locally a compatible complex structure if and only if the
spectrum of the positive Weyl tensor is degenerate
\cite{Nurowski-Przanowski1999}. We refer to \cite{Apo-Gaud1997} for
generalizations of this result in the Riemannian setting and to
\cite{Apostrov1998} for analogous results for arbitrary pseudo-Riemannian
$4$-manifolds.

In this section, we consider a special class of neutral metrics having interesting curvature properties and
provide examples showing that some integrability results in the Riemannian case are not true for neutral
metrics.

A \emph{Walker manifold} is a triple $(M,g,\mathcal{D})$, where $M$ is a smooth manifold, $g$ an indefinite
metric,  and $\mathcal{D}$ a parallel null distribution. The local structure of such manifolds was described by
A.Walker \cite{Walker1950a} in 1950 and we refer to \cite{Derdzinski} for a coordinate-free version of his
theorem. Of special interest are the Walker manifolds admitting a field of null planes of maximum dimension.
Since the dimension of a null plane is not greater than $dim(M)/2$, the lowest possible case is that of neutral
metrics in dimension four admitting a field of parallel null 2-planes.

 Observe that Walker metrics appear as the underlying structure of several specific
pseudo-Riemann\-ian structures such as:

$\bullet$ Hypersymplectic and paraK\"ahler structures \cite{Hit1990},
\cite{kuche}

$\bullet$  Neutral 4-manifolds with parallel real spinor field \cite{Dunajski},
\cite{LawMat2}

$\bullet$ Einstein hypersurfaces in indefinite space forms \cite{Mag1982}

$\bullet$ Indefinite K\"ahler Lie algebras \cite{Ovando2004}

$\bullet$ $2$-step nilpotent Lie groups with degenerate center \cite{Cord-Par}

\vspace{.1in}

 Note also that Walker metrics play a distinguished role in investigating the holonomy of
indefinite metrics (see for example \cite{BBI} and \cite{Ghanam}).

Recall that \cite{Walker1950a}, for every Walker metric $g$ on a $4$-manifold
$M$, there exist local coordinates $(x,y,z,t)$ around any point of $M$ such
that the matrix of $g$ in these coordinates has the following form
\begin{equation}\label{Walker metric}
g_{(x,y,z,t)} =\left(\begin{array}{cccc}
    0&0&1&0\\
    0&0&0&1\\
    1&0&a&c\\
    0&1&c&b\\
\end{array}\right)
\end{equation}
for some smooth functions $a$, $b$ and $c$. In what follows we shall always
assume that the Walker metrics under consideration are given by (\ref{Walker
metric}).

\subsection{Self-dual Walker metrics}
 The components of the curvature tensor of a Walker metric with respect to the
frame $(\partial_x,\partial_y,\partial_z,\partial_t)$ have been computed in
\cite{Ghanam} (see also \cite{Mat,CEM}). They have been used in
\cite{D-RG-RV-L,Dav-Musk2006} to obtain local description of self-dual and
Einstein self-dual Walker metrics.

\begin{thm}\label{Self-duality}
A Walker metric is self-dual if and only if the functions $a, b, c$ have the
form
\[
\begin{array}{lll}
a=x^2yA+x^3B+x^2C+2xyD+xE+yF+G,\\[5pt]
b=xy^2B+y^3A+y^2K+2xyL+xM+yN+P,\\[5pt]
c=x^2yB+xy^2A+x^2L+y^2D+\frac{1}{2}xy(C+K)+xQ+yR+S,
\end{array}
\]
where $A, B, C, etc.$ are arbitrary smooth functions of $(z,t)$.
\end{thm}

\begin{thm}\label {Ein-SD}
A Walker metric is Einstein and self-dual if and only if the functions $a, b,
c$ have the form
\[
\begin{array}{lll}
a=x^2K+xE+yF+G,\\[6pt]
b=y^2K+xM+yN+P,\\[6pt]
c=xyK+xQ+yR+S,
\end{array}
\]
where $K$ is a constant and $E, F, G, etc.$ are arbitrary smooth functions of
$(z,t)$ satisfying the following PDE's:
$$
\begin{array}{lll}
2R_{z}-2F_{t}=FQ+R^2+KG-RE-FN, \\[6pt]
E_{t}+N_{z}-R_{t}-Q_{z}=FM-QR+KS,\\[6pt]
2Q_{t}-2M_{z}=MR+Q^2+KP-EM-QN.
\end{array}
$$
\end{thm}

The next result provides large families of neutral metrics whose hyperbolic twistor spaces are isotropic
K\"ahler but non-K\"ahler (see Theorem \ref{hyperb}\,(iii))

\begin{thm}\label {W-minusB^2|restr}
A Walker metric is self-dual and satisfies the conditions ${\cal
B}^2|\Lambda_{-}=0$ if and only if the functions $a, b, c$ have the form
\[
\begin{array}{lll}
a=x^2K+xE+yF+G,\\[6pt]
b=y^2K+xM+yN+P,\\[6pt]
c=xyK+xQ+yR+S,
\end{array}
\]
where $K, E, F, etc.$  are arbitrary smooth functions of $(z,t)$. In this case
the metric has constant scalar curvature if and only if $K$ is a constant.
\end{thm}

\subsection{Hyperhermitian Walker metrics}

Let $g$ be a Walker metric on ${\Bbb R}^4$ having the form (\ref{Walker
metric}). Then an orthonormal frame of  $T{\Bbb R}^4$ can be specialized by
using the canonical coordinates as follows:
\[
\begin{array}{lll}
{{\bf e}_1=\displaystyle{\frac{1-a}{2}}
\partial_x + \partial_z}, \quad\quad\,\,\, {\bf
e}_2=\displaystyle{\frac{1-b}{2}} \partial_y + \partial_t
-c \partial_x, \\
                                                 \\
{{\bf e}_3=-\displaystyle{\frac{1+a}{2}} \partial_x  + \partial_z, \quad\quad {\bf
e}_4=-\displaystyle{\frac{1+b}{2}}
\partial_y
          +\partial_t - c \partial_x}.
\end{array}
\]
Then the global sections $s_1,s_2,s_3$ of $\Lambda_{-}$ defined by (\ref{eq s-s-basis}) determine via
(\ref{3.1}) a neutral almost hyperhermitian structure $(g,J_1,J_2,J_3)$ on ${\Bbb R}^4$ which is called {\it
proper} in \cite{DDGMMV2}.

Recall that an indefinite almost Hermitian structure $(g,J)$ is said to be
\emph{isotropic K\"ahler} if $\|\nabla J\|^2$ $=$ $0$. Isotropic K\"ahler
structures were first investigated in \cite{Ed-YM2000} in dimension four and
subsequently in \cite{Bla-Dav-Musk2004} in dimension six. It has been shown in
\cite{DDGMMV} that the neutral almost Hermitian structure $(g, J_1)$ on a
Walker 4-manifold is isotropic K\"ahler. Moreover, we have the following
\cite{DDGMMV2}:

\begin{thm}
Any proper almost hyperhermitian   structure $(g,J_1,J_2,J_3)$ on a Walker
$4$-manifold satisfies
 $\|\nabla J_i\|^2= \| d\Omega_i\|^2 = \|\delta\Omega_i\|^2= \| N_i\|^2=0$, where {$N_i$} is the Nijenhuis tensor of
{$J_i$}, {$i=1,2,3$}.
\end{thm}

 Examples of compact isotropic K\"{a}hler structures can be constructed on 4-tori
taking  $a$, $b$ and $c$ in (\ref{Walker metric}) to be periodic functions on $\mathbb{R}^4$. Note also that, in
general, the isotropic K\"ahler structure $(g,J_1)$ is neither K\"ahler nor symplectic.

The next two results have been proved in \cite{DDGMMV2}.

\begin{thm}\label {PC} The structure
$(g,J_1,J_2,J_3)$ is neutral hyperhermitian if and only if the functions $a$,
$b$ and $c$ have the form
\[
\begin{array}{l}
a=x^2K+x P+\xi, \\[0.1in]
b=y^2K+y T+ \eta, \\[0.1in]
c=xyK+\frac{1}{2}x T+\frac{1}{2}y P+ \gamma,
\end{array}
\]
where the capital and Greek letters stand for arbitrary smooth functions of
$(z,t)$.
\end{thm}

\begin{thm}\label{HK}
The structure $(g,J_1,J_2,J_3)$ is neutral hyperk\"ahler if and only if the
functions $a$, $b$ and $c$ do not depend on $x$ and $y$.
\end{thm}

 In particular, the above theorem shows that the metrics considered by Petean \cite{P} are all
neutral hyperk\"ahler and hence self-dual and Ricci-flat. Moreover, this
theorem together with the Petean's classification of neutral  Ricci-flat
K\"ahler surfaces (see Theorem \ref{petean}) leads to the following description
of compact neutral 4-manifolds admitting two parallel, orthogonal and null
vector fields \cite{inprogress}.

\begin{thm}
Let $(M,g)$ be a compact oriented neutral 4-manifold with two parallel ,
orthogonal and null vector fields. Then $M$ is diffeomorphic to a torus or to a
primary Kodaira surface and the metric $g$ is neutral hyperk\"ahler.
\end{thm}

\subsection{Hermitian Walker metrics}
The almost complex structure $J_1$ defined above has been introduced  by
Matsushita \cite{Mat} and
 called the \emph{proper} almost complex structure of the Walker metric $g$. Further we
  denote it by $J$ and note that the structure $(g, J)$ is:

\bigskip
\medskip

{\bf$\bullet$}~ Hermitian  if and only if $a_x-b_x=2c_y, \hspace{.05in}
a_y-b_y=-2c_x.$

\vspace{.2in}

{\bf$\bullet$}~ K\"{a}hler   if and only if $a_x=-b_x=c_y, \hspace{.05in}
a_y=-b_y=-c_x.$

\vspace{.2in}

Moreover, the following results have been proved in \cite{DDGMMV2}.

\begin{thm}
The structure $(g,J)$ is Hermitian and self-dual if and only if
\[
\begin{array}{l}
   a = x^2 K + x y L + x P + y Q + \xi,
   \\[0.1in]
   b = y^2 K - x y L + x S + y T + \eta,
   \\[0.1in]
   c = \frac{1}{2} ( y^2-x^2) L
   + x y K - \frac{1}{2} x (Q-T) + \frac{1}{2} y (P-S) + \gamma,
\end{array}
\]
where all capital and Greek letters are arbitrary smooth functions of $(z,t)$.
\end{thm}

%

\begin{thm}
The structure $(g,J)$ is K\"{a}hler self-dual if and only if
\[
\begin{array}{l}
   a =  x y L + x P + y Q + \xi,
   \\[0.1in]
   b = - x y L - x P - y Q + \eta,
   \\[0.1in]
   c = \frac{1}{2} (y^2 - x^2) L     - x Q + y P + \gamma.
\end{array}
\]

\end{thm}

\begin{thm}\label{Kahler-Einstein}
The structure $(g,J)$ is K\"{a}hler-Einstein  if and only if

\[
\begin{array}{lll}
   a=\kappa(x^2- y^2) + x P +y Q+ \xi,
   \\[6pt]
   b=\kappa(y^2-x^2) - x P -y Q   -\xi + \frac{1}{\kappa}( P_{z}-Q_{t}),
   \\[6pt]
   c=  2\kappa xy  - x Q  + y P + \gamma,
\end{array}
\]

or

\[
\begin{array}{lll}
a= x P+y Q+ \xi,
\\[6pt]
b= -x P-y Q+ \eta,
\\[6pt]
c= -x Q +y P+\gamma,
\end{array}
\]
where in the last case
   $ P_z=Q_t.$
Here $\kappa$ is a non-zero constant and all capital and Greek letters are
smooth functions of $(z,t)$. In the first case $\tau=8\kappa$, in the second
one $\tau=0$.
\end{thm}

The above theorem has been used in \cite{DDGMMV2} to construct examples of
neutral almost Hermitian Einstein metrics showing that an integrability result
of Kirchberg \cite{Kirch2006} in the positive definite case is not valid in the
neutral setting.

\subsection{Almost K\"ahler Walker metrics}

Although the Goldberg conjecture \cite{Go} mentioned above is of global nature, it is already known that some
additional curvature conditions suffice to show the integrability of the almost complex structure at the local
level. For instance, in dimension four, Einstein almost K\"{a}hler metrics which are $\ast$-Einstein are
necessarily K\"{a}hler \cite{Ogr-Seki1998}. (The $\ast$-Einstein condition can be replaced by the second Gray
curvature identity or the anti-self-duality condition and the integrability still follows \cite{Apo-Dra2003}).
It is also well-known \cite{Arm2002, Ogr-Seki1998, Ols1978}  that any almost K\"ahler metric of constant
sectional curvature is K\"ahler (and flat).

 Our purpose here is to exhibit large families of neutral strictly
almost K\"ahler Einstein metrics showing that the results mentioned above are
not true for neutral metrics. Their construction is based on the following
result \cite{DDGMMV}:

\begin{thm}\label {Ein-AK}
A proper structure $(g,J)$ is strictly almost K\"ahler Einstein if and only if
the functions $a$, $b$ and $c$ have the form
\begin{align}
& a=xP+yQ+\xi,  \nonumber \\[6pt]
& b=-xP-yQ+\eta, \nonumber \\[6pt]
& c=xU+yV+\gamma ,  \nonumber
\end{align}
where all capital and Greek letters are smooth functions of $(z,t)$ satisfying
the following PDE's:
\[
\begin{array}{lll}
2(V_{z}-Q_{t})=V^2-VP+Q^2+UQ,\\[6pt]
2(P_{z}+U_{t})=P^2-VP+U^2+UQ,\\[6pt]
Q_{z}+ U_{z}-P_{t}+V_{t}=PQ+UV,
\end{array}
\]
and $(V-P)^2+(U+Q)^2\not\equiv 0$.
\end{thm}

An interesting consequence of this result (compare with the Riemannian case) is
the folowing:

\begin{cor}\label {AK}
Any proper strictly almost K\"ahler Einstein structure on a Walker 4-manifold
is self-dual, Ricci flat and $\ast$-Ricci flat.
\end{cor}

{\bf Example} For any constants $p, q, r$ with $p^2+q^2\neq 0$, set
\begin{align}
a=-\frac{2(px-qy)}{pz+qt+r},\quad b=-\frac{2(px+qy)}{pz+qt+r} \nonumber.
\end{align}
 It has been shown in \cite{DDGMMV} that the proper almost Hermitian
structure $(g,J)$ defined by the functions $a, b$ and $c=0$ is strictly almost
K\"ahler and the metric $g$ is flat.

Finally, let us note that we do not know of compact examples of neutral strictly almost K\"ahler-Einstein
4-manifolds. However, an indefinite Ricci-flat strictly almost Kahler metric on 8-dimensional torus has been
recently reported in \cite{MHL}.

\section{Neutral anti-self-dual 4-manifolds and projective \\ \centerline {structures}}

  In  \cite{DW1} Dunajski and West gave a local classification of neutral anti-self-dual
  4-manifolds admitting a non-trivial null conformal Killing vector
field. It is based on the observation that the self-dual and anti-self-dual
null plane distributions associated to such a vector field  are integrable.
Then their main result is that the two dimensional leaf space of the foliation
determined by the anti-self-dual plane distribution is equipped with a
canonical projective structure, i.e. a class of torsion-free connections with
same unparametrized geodesics. Conversely, any projective structure on a
2-dimensional surface gives rise to a neutral anti-self-dual 4-manifold with
null conformal Killing vector field. The explicit form of this correspondence
is given by the following:

\begin{thm} Let $(M, [g], K)$ be a neutral anti-self-dual conformal
structure with a
 null conformal Killing vector field $K$.
  Then around any point where $K\neq 0$ there exist local coordinates
  $(\phi,x,y,z)$ and $g\in [g]$ such that $K=\frac{\partial}{\partial\phi}$
   and $g$ has one of the following two forms, according to whether the
   twist ${\Bbb K}\wedge d{\Bbb K}$ vanishes or not (here ${\Bbb K}$ is the
   one-form defined by ${\Bbb K}:=g(K,.)$):

 \item[$(i)$] If ${\Bbb K}\wedge d{\Bbb K}=0$, then
\[
\begin{array}{lll}
 g&=&(d\phi+(zA_3-Q)dy)(dy-\beta dx)\\
 & &-(dz-(z(-\beta_y+A_1+\beta A_2+\beta^2A_3))dx-(z(A_2+2\beta A_3)+P)dy)dx,
\end{array}
\]
  where $A_1,A_2,A_3,\beta,P,Q $ are arbitrary functions of $(x,y)$.

 \item[$(i)$] If ${\Bbb K}\wedge d{\Bbb K}\neq 0$, then
 $$
g=(d\phi+A_3\frac{\partial G}{\partial z}dy+(A_2\frac{\partial G}{\partial
z}+2A_3(z\frac{\partial G}{\partial z}-G) -\frac{\partial^2 G}{\partial
z\partial y})dx)(dy-zdx)
$$
$$
-\frac{\partial^2 G}{\partial z^2}dx(dz-(A_0+zA_1+z^2A_2+z^3A_3)dx),
$$
where $A_0,A_1,A_2,A_3$ are arbitrary functions of $(x,y)$, and $G$ is a
function of $(x,y,z)$ satisfying the following PDE:
 $$(\frac{\partial}{\partial z}+\frac{\partial}{\partial y}+(A_0+zA_1+z^2A_2+z^3A_3)\frac{\partial}{\partial z})\frac{\partial^2 G}{\partial
 z^2}=0.$$

 \end{thm}

 The above two forms of the metric $g$ are related to the two dimensional
projective structures in the following way. The two geodesic's equations for
any connection in such a projective class can be reduced, by eliminating the
affine parameter, to one second order ODE of the following form:
$$
\frac{d^2y}{dx^2}=A_3(x,y)(\frac{dy}{dx})^3+A_2(x,y)(\frac{dy}{dx})^2+A_1(x,y)(\frac{dy}{dx})
+A_0(x,y).$$

 The functions $A_i$ can be expressed in terms of the connection
coefficients and do not depend on the particular choice of the connection.
These are the $A_i$'s which appear in the theorem above, where in the first
case we have $A_0=\beta_x+\beta\beta_y-\beta A_1-\beta^2A_2 -\beta^3A_3$.

The relation between the anti-self-dual neutral 4-manifolds with a non-trivial null conformal Killing vector
field and the two dimensional projective structures found by Dunajski and West has been generalized by
Calderbank \cite{Cald} for neutral anti-self-dual 4-manifolds admitting an integrable anti-self-dual null
distribution instead of a null conformal Killing vector field. We should also mention the recent work by Bryant,
Dunajski and Eastwood \cite{BDE} on two dimensional metrisable projective structures which suggests the
existence of a conformal invariant on neutral anti-self-dual manifolds, which vanishes when the conformal class
contains a neutral K\"ahler metric with conformal null symmetry.

 \vspace{0.5cm}

\section{\bf Twistor spaces}

One can construct a twistor space of an oriented 4-manifold with a neutral metric as in the Riemannian case
\cite{BDM} (see also \cite{AS}, where the pseudo-sphere $SO(3,2)/SO(2,2)$ with the standard invariant neutral
metric is discussed). These twistor spaces are called hyperbolic in \cite{BDM} since their fibre is a two
sheeted hyperboloid.

In this section we discuss the hyperbolic twistor spaces as well as their
relation to the LeBrun-Mason twistor spaces of neutral self-dual manifolds
whose fibre is ${\Bbb C}{\Bbb P}^1$ \cite{LB-M}.

\subsection{Hyperbolic twistor spaces}

  Let $M$ be an oriented 4-manifold with a
neutral metric $g$. Then the hyperbolic twistor space ${\cal Z}$ of $M$ is
defined to be the bundle on $M$ consisting of all unit bivectors in
$\Lambda_-$. It can be identified via (\ref{3.1}) with the space of all complex
structures on the tangent spaces of $M$ compatible with its metric and
orientation. Note that the fiber of the twistor bundle $\pi: {\cal
Z}\rightarrow M$ is the two sheeted hyperboloid
$$ H=\{(y_1,y_2,y_3)\in {\Bbb R}^3:\,-y_1^2+y_2^2+y_3^2=-1\}$$
and we can think of the $y_1>0$ branch as one of the standard models of the
hyperbolic plane. Further, we shall consider the hyperboloid $H$ with the
complex structure $S$ determined by the restriction to $H$ of the metric
$-dy_1^2+dy_2^2+dy_3^2$, i.e. $SV = y\times V$ for $V\in T_yH$, where $\times$
is the vector cross product on ${\Bbb R}^3$ defined by means of the
paraquaternionic algebra.

 The Levi-Civita connection of $g$ gives rise to a splitting $T{\cal Z}={\cal
 H}\oplus {\cal V}$ of the tangent bundle of ${\cal Z}$ into horizontal and
 vertical components and we consider the vertical space ${\cal V}_{\sigma}$ at
 $\sigma\in {\cal Z}$ as the orthogonal complement of $\sigma$ in
 $\Lambda_-$. As in the Riemannian case \cite{AHS, ES}, we define two almost complex structures
  ${\cal J}_1$ and ${\cal J}_2$ on ${\cal Z}$ by
  $$ {\cal J}_nV= (-1)^{n-1}\sigma\times V ,\hspace{.1in} V\in {\cal V}_{\sigma},$$
  $$\pi_*({\cal J}_n A) = J_{\sigma}(\pi_*A) , \hspace{.1in} A\in {\cal H}_{\sigma},$$
where $J_{\sigma}$ is the complex structure on $T_pM$, $p\in\pi(\sigma)$,
  defined by (\ref{3.1}).

    The metric $g$ on the bundle $\pi:\Lambda^2
\longrightarrow M$ induced by the metric of $M$ is negative definite on the fibres of ${\cal Z}$ and we adopt
the metric $\langle\,,\rangle=-g$ on $\Lambda^2$. Setting $h_t=\pi^*g +t\langle\,,\rangle$ for any real $t\neq
0$, we get a $1$-parameter family of pseudo-Riemannian metrics on ${\cal Z}$ compatible with the almost complex
structures ${\cal J}_1$ and ${\cal J}_2$. The almost Hermitian structures $(h_t,{\cal J}_k), t\neq 0, k=1,2$,
have been studied in \cite{Bla-Dav-Musk2004,BDM} where the following results have been proved:

\begin{thm}\label{hyperb}
On the hyperbolic twistor space ${\cal Z}$ of an oriented $4$-manifold $M$ with
a neutral metric $g$ we have the following:
\begin{enumerate}
\item[$(i)$] The almost complex structure ${\cal J}_1$ is integrable if and
only if the metric $g$ is self-dual. The almost Hermitian structure $({\cal
J}_1,h_t)$ is indefinite K\"ahler if and only if $g$ is an Einstein and
self-dual metric with constant scalar curvature $\tau= -12/t$. \item[$(ii$)]
The almost complex structure ${\cal J}_2$ is never integrable. The almost
Hermitian structure $({\cal J}_2,h_t)$ is indefinite almost K\"ahler (resp.
nearly K\"ahler) if and only if the metric $g$ is Einstein, self-dual and
$\tau=12/t$ (resp.  $\tau=-6/t$). \item[$(iii$)] The almost Hermitian structure
$(h_t,{\cal J}_k)$ is isotropic K\"ahler if and only if $k=1$, $g$ is a
self-dual metric with constant scalar curvature $\tau=-12/t$ and ${\cal
B}^2{|\Lambda_-}=0$.
\end{enumerate}
\end{thm}

 Note  that the values of the scalar curvature appearing in Theorem \ref{hyperb} for $t>0$ are the
negatives of what one has for the usual twistor space \cite{DMG}. This sign
change is due to our choice of metric on $\Lambda^2$.

Statement $(iii)$ of the above theorem motivates the study of self-dual neutral metrics of constant scalar
curvature and two-step nilpotent Ricci operator ${\cal B}$. Note that it follows from Theorem
\ref{W-minusB^2|restr} that this class of neutral metrics is strictly larger than the class of neutral self-dual
Einstein metrics. The next example shows that it contains also neutral conformally flat metrics of non-constant
sectional curvature.

\medskip

 \noindent {\bf  Example}
Let $G$ be a Lie group whose Lie algebra has a basis $\{E_1, E_2, E_3, E_4\}$
such that
$$[E_1,E_2]= E_2,\quad
[E_1,E_3]=-E_2+ 3 E_3,\quad [E_1,E_4]=2 E_4,$$
$$[E_2,E_3]=[E_2,E_4]= [E_3,E_4]=0.$$  Define a left-invariant neutral metric $g$ on $G$ in
terms of the dual basis $\{E^i\}$ by
$$g=E^1\otimes E^1+E^2\otimes E^2-E^3\otimes E^3-E^3\otimes E^3.$$
Then it is straightforward to compute the curvature operator ${\cal R}$ of $g$
and to see that $\tau = -12$, ${\cal W}=0$, ${\cal B}^2=0$ but ${\cal B}\neq
0$. More examples of such metrics have been constructed in \cite{Dav-Musk2006}.

\medskip

Next we consider the problem for local and global existence of holomorphic
functions on hyperbolic twistor spaces. On the classical twistor space over an
oriented Riemannian 4-manifold with either almost
  complex structure, there are no global non-constant
holomorphic functions, even
  when the base manifold is
non-compact \cite{DM, DMG}.  However for the hyperbolic twistor spaces there is
  considerable difference from the classical case as we
shall see.

Recall that a $C^\infty$ complex-valued function on an almost complex manifold
is said to be {\it holomorphic} if its differential is complex-linear with
respect to the almost complex structure. For any $n=0,1,2,3$, let ${\cal
F}_n({\cal J}_i)$ denote the (possibly empty) set of points $\sigma\in {\cal
Z}$ such that $n$ is the maximal number of local ${\cal J}_i$-holomorphic
functions with ${\Bbb C}$-linearly independent differentials at $\sigma$.  In
\cite{DM} it is shown that  for the classical twistor space ${\cal Z}$ of an
oriented Riemannian 4-manifold $(M,g)$ we have ${\cal Z}={\cal F}_0({\cal
J}_1)\cup{\cal F}_3({\cal J}_1) ={\cal F}_0({\cal J}_2)\cup{\cal F}_1({\cal
J}_2)$; moreover
$${\cal F}_3({\cal J}_1)=\pi^{-1}(Int\{p\in M:\, ({\cal
W}_{-})_p=0\}),$$
$${\cal F}_1({\cal J}_2)=\pi^{-1}(Int\{p\in M:\, {\cal
R}_p=({\cal W}_{+})_p\}).$$ The same arguments give this result for the
hyperbolic twistor space as well.

   Given a neutral almost hyperhermitian
4-manifold $(M,g,J_1,J_2,J_3)$, denote by $\pi:{\cal Z}\to M$ the hyperbolic
twistor space of $(M,g)$. Then the $2$-vectors corresponding to $J_1,J_2,J_3$
via (\ref{3.1}) form a global frame of $\Lambda_-$ and we have a natural
projection $p:{\cal Z}\to H$ defined by $p(\sigma)=(y_1,y_2,y_3)$, where
$J_{\sigma}=y_1J_1(x)+y_2J_2(x)+y_3J_3(x)$, $x=\pi(\sigma)$. Thus ${\cal Z}$ is
diffeomorphic to $M\times H$ by the map $\sigma\to (\pi(\sigma),p(\sigma))$ and
it is obvious that $p$ maps any fibre of ${\cal Z}$ biholomorphically on $H$
with respect to ${\cal J}_1$ and $S$. The conditions on $(M,g)$ ensuring that
the natural projection $p:{\cal Z}\to H$ is ${\cal J}_1$-holomorphic (resp.
${\cal J}_2$-anti-holomorphic) are the following \cite{BDM}:

\begin{thm}
    Let $M$ be a neutral almost hyperhermitian
4-manifold and ${\cal Z}$ its hyperbolic twistor space. Then the natural
projection $p:{\cal Z}\to H$ is ${\cal J}_1$-holomorphic (resp. ${\cal
J}_2$-anti-holomorphic) if and only if $M$ is neutral hyperhermitian (resp.
neutral hyperk\"ahler).
\end{thm}

Furthermore we have the following result \cite{BDM} in the compact case:

\begin{thm}
   Let $M$ be a compact neutral hyperhermitian manifold
with hyperbolic twistor space ${\cal Z}$ and let $p:{\cal Z}\to H$ be the natural projection. Then any ${\cal
J}_1$-holomorphic function $f$ on ${\cal Z}$ has the form $f=g\circ p$, where $g$ is a holomorphic function on
$H$. If $M$ is neutral hyperk\"ahler, any ${\cal J}_2$-holomorphic function $f$ on ${\cal Z}$ has the form
$f=g\circ p$, where $g$ is an anti-holomorphic function on $H$.
\end{thm}

Now we show that in the non-compact case the situation changes drastically. Recall that the Petean metrics on
${\Bbb R}^4$ have the form:
$$g=f(dx_1\otimes dx_1+dx_2\otimes dx_2)+ dx_1\otimes
dx_3+ dx_3\otimes dx_1 + dx_2\otimes dx_4 + dx_4\otimes dx_2,$$ where $(x_1,x_2,x_3,x_4)$ are the standard
coordinates on ${\Bbb R}^4$ and $f$ is a smooth positive function depending on $x_1$ and $x_2$ only. According
to Theorem \ref{HK} the metrics $g$ are neutral hyperk\"ahler, i.e. self-dual and Ricci-flat. Hence the almost
complex structure ${\cal J}_1$ on the hyperbolic twistor space ${\cal Z}$ of $({\Bbb R}^4,g)$ is integrable.
Moreover, we have the following result \cite{BDM}, which shows that there can be an abundance of global
holomorphic functions on a hyperbolic twistor space.
\begin{thm}
The hyperbolic twistor space $({\cal Z},{\cal J}_1)$ of $({\Bbb R}^4,g)$ is
biholomorphic to ${\Bbb C}^2\times H$.
\end{thm}

This result suggests the following problem: {\it Characterize the non-compact neutral hyperk\"ahler 4-manifolds
whose hyperbolic twistor spaces $({\cal Z},{\cal J}_1)$ are Stein manifolds.}

\subsection{LeBrun-Mason twistor spaces}

As we noted in the previous subsection, the hyperbolic twistor spaces are non-compact, hence the problem of
their compactification arises. Since the fibre is a two sheeted hyperboloid, it is natural to compactify it by
means of the sphere ${\Bbb C}{\Bbb P}^1$.
 To put this simple observation in a global framework, we shall follow a recent
 twistor construction  by LeBrun and Mason \cite{LB-M}.

 Denote by $g$ the complex-bilinear extension of $g$ to the complexification $T^{\Bbb C}M$ of the tangent bundle
 $TM$. Let $\Lambda_{-}^{\Bbb C}=\Lambda_{-}\otimes{\Bbb C}$ be the complexification of $\Lambda_{-}$. Set
$$
{\cal Z}=\{[\varphi]\in{\Bbb P}(\Lambda_{-}^{\Bbb C}):\, g(\varphi,\varphi)=0\}
$$
 and let ${\cal P}:{\cal Z}\to M$ be the natural projection. Then ${\cal P}:{\cal Z}\to M$ is a ${\Bbb C}{\Bbb
P}^1$-bundle. Indeed, if $\{e_1, e_2, e_3, e_4\}$ is an oriented orthonormal frame of vector fields on an open
set $U$, denote by $s_1, s_2, s_3$ the frame of $\Lambda_{-}$ defined by (\ref{eq s-s-basis}). Then
$$
\Psi:([\zeta_1,\zeta_2],x)\to
([(\zeta_1^2+\zeta_2^2)s_1(x)+(\zeta_1^2-\zeta_2^2)s_2(x)-2\zeta_1\zeta_2
s_3(x)],x),\> x\in U,
$$
is a biholomorphic map of ${\Bbb C}{\Bbb P}^1\times U$ onto ${\cal P}^{-1}(U)$.

The bundle
$$
F=\{[\varphi]\in {\Bbb P}(\Lambda_{-}):\, g(\varphi,\varphi)=0\}
$$
is embedded in ${\cal Z}$ in an obvious way. Each fibre $F_x$ of $F$ is sent by $\Psi^{-1}$ onto ${\Bbb R}{\Bbb
P}^1\times \{x\}$. Note also that ${\Bbb R}{\Bbb P}^1$ goes to the "Greenwich" meridian $S^1:\>
\xi_1^2+\xi_2^2=1,\>\xi_3=0$ under the standard identification ${\Bbb C}{\Bbb P}^1\cong S^2$. The manifold $S^2
- S^1$ is biholomorphic to the disjoint union of two copies of the unit disk and the latter manifold is
diffeomorphic by the stereographic projection to the two sheeted hyperboloid $\xi_1^2-\xi_2^2-\xi_3^2=1$. Thus
${\cal Z}_x - F_x$ can be identified with the fibre of the hyperbolic twistor space of $(M,g)$. This can be also
seen in the following way. If $\Pi\in {\cal Z}_x - F_x$, then $\Pi\cap\overline\Pi=0$, hence $T^{\Bbb
C}_xM=\Pi\oplus\overline\Pi$. Denote by $J_x=J_x(\Pi)$ the (unique) complex structure on the vector space $T_xM$
for which $\Pi$ is the space of the $(1,0)$-vectors. The fact that the space $\Pi$ is isotropic is equivalent to
the compatibility of $J_x$ with the metric $g$. If $\{e_1,e_2=J_xe_1,e_3,e_4=J_xe_3\}$ is a basis of $T_xM$ with
$||e_1||^2=1$, $||e_3||^2=-1$, then $\Pi$ is represented by the bivector $(e_1-ie_2)\wedge (e_3-ie_4)=(e_1\wedge
e_3-e_2\wedge e_4)-i(e_1\wedge e_4+e_2\wedge e_3)$. This vector lies in $\Lambda^{\Bbb C}_{-}$ since
$\Pi\in{\Bbb P}(\Lambda^{\Bbb C}_{-})$, therefore the above basis induces the orientation of $T_xM$, which means
that $J_x$ is compatible with this orientation. Thus the assignment ${\cal Z}_x - F_x\ni\Pi\to J_x(\Pi)$
identifies  ${\cal Z}-F$ with the hyperbolic twistor space of $(M,g)$.

 We have assumed that the manifold $M$ is oriented, i.e. its structure group is reduced to $SO(2,2)$. Now
suppose that the structure group of $M$ can be reduced to the identity
component $SO_{+}(2,2)$ of $SO(2,2)$. In this case, following \cite{LB-M}, we
shall say that $M$ is space-time orientable. This condition is equivalent to
the existence of two $g$-orthogonal subbundles $T_{\pm}$ of $TM$ such that
$TM=T_{+}\oplus T_{-}$ and the restriction of $g$ to $T_{+}$, resp. $T_{-}$, is
positive, resp., negative definite. To see this choose a Riemannian metric $h$
on $M$, then diagonalize $g$ with respect to $h$.

 Note that if the bundle $T_{\pm}$ is oriented, we can define a unique complex structure ${\cal J}_{\pm}$ on $T_{\pm}$
compatible with the metric and the orientation since $rank\> T_{\pm}=2$ and the
restriction of $g$ on $T_{\pm}$ is a definite bilinear form: if $e_1, e_2$ is
an oriented orthonormal basis of $T_{\pm}$, ${\cal J}_{\pm}e_1=e_2$.

By definition, to fix a space-time orientation of $M$ means to choose
orientations on the bundles $T_{\pm}$. In this case $TM$ is considered with the
orientation obtained via the decomposition $TM=T_{+}\oplus T_{-}$.

  Suppose that $(M,g)$ is space-time oriented and denote by ${\cal J}$ the almost complex structure on $M$, which
coincides with ${\cal J}_{\pm}$ on $T_{\pm}$. Let $\{e_1, e_2={\cal J}e_1\}$ and $\{e_3, e_4={\cal J}e_3\}$ be
oriented orthonormal bases of $(T_{+})_x$ and $(T_{-})_x$, respectively, $x\in M$. Then the $2$-plane of
$(1,0)$-vectors of ${\cal J}_x$ is represented by the bivector $(e_1-ie_2)\wedge (e_3-ie_4)=s_2-is_3$, while the
2-plane corresponding to the conjugate complex structure $-{\cal J}_x$ is represented by $s_2+is_3$. The images
of $[s_2-is_3]$ and $[s_2+is_3]$ under $\Psi^{-1}$ are the points $[1,i]$ and $[1,-i]$ in ${\Bbb C}{\Bbb P}^1$,
which go to $(0,1,0)$ and $(0,-1,0)$ under the standard identification ${\Bbb C}{\Bbb P}^1\cong S^2$. The latter
points determine the two connected components of $S^2 - S^1$. Thus, if $(M,g)$ is space-time oriented, the space
${\cal Z} - F$ has two connected components - one of them contains the almost complex structure ${\cal J}$ and
the other one contains the conjugate almost complex structure $-{\cal J}$.  Let $U$ be the component of ${\cal
Z} - F$ determined by ${\cal J}$. Then the closure ${\cal Z}_{+}=U\cup F$ of $U$ in ${\cal Z}$ is a compact
$6$-dimensional manifold with boundary.

 Now suppose that the image of every maximally extended null geodesic in $M$ is an embedded circle in $M$. A
manifold with this property is called Zollfrei and we refer to \cite{LB-M} for more information about these
manifolds. It is proved in \cite{LB-M} that if $(M,g)$ is a space-time oriented self-dual Zollfrei $4$-manifold,
the space $F$ is a trivial $2$-sphere bundle $q:F\to P$ over a manifold $P$ diffeomeorphic to ${\Bbb R}{\Bbb
P}^3$. Let us endow the disjoint union $Z=U\coprod P$ with the quotient topology induced by the map ${\cal
Z}_{+}\to Z$, which is the identity map on $U$ and the map $q$ on $F$. Then $Z$ is a compact topological
$6$-manifold. One of the key results in \cite{LB-M} is that $Z$ admits a complex structure, which coincides with
the Atiyah-Hitchin-Singer structure ${\cal J}_1$ on $U$, identified with a connected component of the hyperbolic
twistor space of $M$. The manifold $Z$ is called the twistor space of $(M,g)$.

 If $(M,g)$ is an oriented self-dual Zollfrei manifold, which is not space-time orientable, it possesses a double
cover $(\widetilde M,\widetilde g)$, which is space-time orientable, self-dual and Zollfrei. The twistor space
of $(\widetilde M,\widetilde g)$ is called the twistor space of $(M,g)$.

  It is proved in \cite{LB-M} that if $(M,g)$ is a self-dual Zollfrei $4$-manifold and if $M$ is space-time
orientable, then it is homeomorphic to $S^2\times S^2$; if $M$ is not
space-time orientable, it is homeomrphic to the quadric
$$
{\Bbb M}^{2,2}=\{[x,y]\in{\Bbb R}{\Bbb P}^5:\, |x|^2-|y|^2=0\},
$$
which is the quotient of $S^2\times S^2$ by the ${\Bbb Z}_2$-action generated
by the map $(x,y)\to (-x,-y)$. It is also shown in \cite{LB-M} that the twistor
space of either $S^2\times S^2$ or ${\Bbb M}^{2,2}$ with the natural neutral
metric is biholomorphic to ${\Bbb C}{\Bbb P}^3$ in such a way that $P\subset Z$
becomes the standard ${\Bbb R}{\Bbb P}^3\subset {\Bbb C}{\Bbb P}^3$. Hence the
hyperbolic twistor space of $S^2\times S^2$ is ${\Bbb C}{\Bbb P}^3 - {\Bbb
R}{\Bbb P}^3$. It is also a result of \cite{LB-M} that the twistor space of any
self-dual Zollfrei manifold is diffeomorphic to ${\Bbb C}{\Bbb P}^3$ in such a
way that the Chern classes of $Z$ are sent to the Chern classes of ${\Bbb
C}{\Bbb P}^3$.

  Finally, let us note that the LeBrun-Mason twistor construction can be reversed to produce
self-dual metrics of neutral signature. We refer to \cite{LB-M} for details.

 \vskip 20pt

\noindent Johann Davidov and Gueo Grantcharov

\noindent Institute of Mathematics and Informatics

\noindent Bulgarian Academy of Sciences

\noindent 1113 Sofia, Bulgaria

\noindent and

\noindent Department of Mathematics

\noindent Florida International University

\noindent Miami, FL 33199

\noindent jtd@math.bas.bg,  grantchg@fiu.edu

\medskip

\noindent Oleg Mushkarov

\noindent Institute of Mathematics and Informatics

\noindent Bulgarian Academy of Sciences

\noindent 1113 Sofia, Bulgaria

\noindent muskarov@math.bas.bg


\begin{thebibliography}{99}

\bibitem{AS} R. Albuquerque, I. M. Salavessa,  On the twistor space of pseudo-spheres,
{\it Diff. Geom. Appl.} {\bf 25} (2007), 207-219.





\bibitem{Apo-Gaud1997}
V. Apostorov and P. Gauduchon, {The Riemann Goldberg-Sacks theorems}, {\it
Internat. J. Math.} {\bf 8} (1997), 421--439.


\bibitem{Apostrov1998}
V.~Apostorov, Generalized Goldberg-Sacks theorems for pseudo-Riemannian
four-manifolds, {\it J. Geom. Phys.} {\bf 27} (1998), 185--198.



\bibitem{Apo-Dra2003}
V. Apostorov, T. Draghici, The curvature and the integrability of almost
K\"ahler manifold: a survey. \textit{Symplectic and contact topology:
interactions and perspectives. } Toronto, ON/Montreal, QC, 2001), 25--53,
Fields Inst. Commun., {\bf 35} Amer. Math. Soc., Providence, RI, 2003.


\bibitem{Arm2002} J. Armstrong, An ansatz for Almost-K\"ahler, Einstein $4$-manifolds,
{\it J. Reine Angew. Math.} {\bf 542} (2002), 53--84.


\bibitem{AHS}
M. Atiyah, N. Hitchin, I. Singer, Self-duality  in four-dimensional Riemannian
geometry, {\it Proc. Roy. Soc. London}  Ser. A {\bf 362} (1978), 425-461.

\bibitem{BPV} W. Barth, K. Hulek, C. Peters, A. Van de Ven, { Compact complex surfaces},
2nd eddition, Springer, 2004.

\bibitem{bes}
A. Besse, { Einstein Manifolds}, Springer, Berlin, 1987.


\bibitem{BBI}
L. B\'{e}rard Bergery, A. Ikemakhen, Sur l'holonomie des vari\'{e}t\'{e}s
pseudo-Riemanniennes de signature $(n,n)$, {\it Bull. Soc. Math. France} {\bf
125} (1997), 93--114.



\bibitem{Bla-Dav-Musk2004}
D. Blair, J. Davidov, O. Mu\v skarov, Isotropic K\"ahler hyperbolic twistor
spaces, \textit{J. Geom. Phys.} {\bf 52} (2004), 74--88.




\bibitem{BDM}
D. Blair, J. Davidov, O. Mushkarov, Hyperbolic twistor spaces, {\it Rocky
Mountain J.Math.} {\bf 35} (2005), 1437-1465.



\bibitem{BDE} R. Bryant, M. Dunajsky, M. Eastwood, Metrizability of
two-dimensional projective structures, preprint arXiv: 0801.0300.

\bibitem{CEM} M. Chiachi, E.~Garc\'{\i}a-R\'{\i}o, Y. Matsushita, Curvature
properties of four-dimensional Walker metrics, {\it Class. Quantum
 Grav.} {\bf 22} (2005), 559-577.


\bibitem{Cald} D. Calderbank, Self-dual 4-manifolds, projective structures and
the Dunajsky-West construction, preprint arXiv: math/0606754.


\bibitem{Cord-Par}
L.~A.~Cordero, P.~F.~Parker, Pseudoriemannian 2-step nilpotent Lie groups,
preprint, arXiv: math.DG/9906188.


\bibitem{Uhl} B. Dai, C.-C. Terng, K. Uhlenbeck, On space-time monopole
equation, preprint, arXiv: math.DG/0602607.



\bibitem{DDGMMV}
J.~Davidov, J. C.~D\'{\i}az-Ramos, E.~Garc\'{\i}a-R\'{\i}o, Y.~Matsushita,
O.~Mu\v skarov, R.~V\'{a}zquez-Lorenzo, Almost K\"{a}hler Walker $4$-manifolds,
{\it J. Geom. Phys.} {\bf 57} (2007), 1075--1088.

\bibitem{DDGMMV2}
J.~Davidov, J. C.~D\'{\i}az-Ramos, E.~Garc\'{\i}a-R\'{\i}o, Y.~Matsushita, O.~Mu\v skarov,
R.~V\'{a}zquez-Lorenzo, Hermitian Walker $4$-manifolds {\it J. Geom. Phys.}, {\bf 58} (2008), 307-323.


\bibitem{DM}
J. Davidov,  O. Mu\u skarov,  Existence of holomorphic functions on twistor
spaces, {\it  Bull. Soc. Math. Belgique} {\bf 40} Ser. B (1989), 131-151.


\bibitem{DMG}
J. Davidov, O. Mu\u skarov, G. Grantcharov,  Almost complex structures on
twistor spaces, {\it Almost Complex Structures}, World Scientific, Singapore,
(1994), 113-149.

\bibitem{inprogress}
J. Davidov, G. Grantcharov, O. Mushkarov, Work in progress.

\bibitem{Dav-Musk2006}
J.~Davidov, O.~Mu\v skarov, Self-dual Walker metrics with a two-step nilpotent
Ricci operator, {\it  J. Geom. Phys.} {\bf 57} (2006), 157--165.

\bibitem{Derdzinski} A. Derdzinski, W. Rotter, Walker's theorem without coordinates,
{\it J. Math. Phys.} {\bf 47} (2006), no. 6, 062504, 8 pp.

\bibitem{D-RG-RV-L} J.~Carlos D\'{\i}az-Ramos, Eduardo Garc\'{\i}a-R\'{\i}o, Ram\'{\o}n
V\'azquez-Lorenzo,  {Four-dimensional Osserman metrics with non-diagonalizable Jacobi operator}, {\it J. Geom.
Anal.} {\bf 16} (2006), 39-52.

\bibitem{Dunajski} M. Dunajski, Anti-self-dual four-manifolds with a parallel
real spinor {\it Proc. R. Soc. London A} {\bf 458} (2002) 1205-1222.


\bibitem{DW1} M. Dunajski, S.West,  Anti-self-dual conformal
 structures with null Killing vectors from projective structures.
 {\it Comm. Math. Phys.} {\bf 272} (2007), 85--118.

\bibitem{DW} M. Dunajsky, S. West,  Anti-Self-Dual Conformal Structures in Neutral
Signature, to appear in the special volume {\it `Recent developments in
pseudo-Riemannian Geometry'}, ESI-Series on Mathematics and Physics, preprint,
arXiv: math/0610280.

\bibitem{ES}
J. Eells, S. Salamon, Twistorial construction of harmonic maps of surfaces into
four-manifolds, {\it Ann. Scuola Norm. Sup.  Pisa}  Cl.\,Sci. {\bf 12} (1985),
589-640.



\bibitem{Ed-YM2000}
E.~Garc\'{\i}a-R\'{\i}o, Y.~Matsushita, {Isotropic K\"ahler structures on Engel
4-manifolds}, \textit{J. Geom. Phys.} {\bf 33} (2000), 288--294.

\bibitem{Ghanam} R. Ghanam, G. Tompson, The holonomy Lie algebras of neutral
metrics in dimension four, {\it J. Math. Phys.} {\bf 42} (2001), 2266-2285.


\bibitem{Go}
S. Goldberg, Integrability of almost K\"{a}hler manifolds, \emph{Proc. Amer.
Math. Soc.} \textbf{21} (1969), 96--100.


\bibitem{MHL}
S. Haze, Y.~Matsushita, P. Law, Almost K\"ahler-Einstein structures on
8-dimensional Walker manifolds, {\it Monatsh. Math.} {\bf 150} (2007), 41-48.



\bibitem{Hit1990}
N.~Hitchin, {Hypersymplectic quotients}, \textit{Acta Acad. Sci. Tauriensis}
{\bf 124} (1990), 169--180.





\bibitem{kuche} S. Ivanov, S. Zamkovoy, Parahermitian and paraquaternionic
manifolds, {\it Diff. Geom. Appl.} {\bf 23} (2005), 205--234.


\bibitem{Ka}
H. Kamada, Neutral hyperk\"ahler structures on primary Kodaira surfaces, {\it
Tsukuba J. Math.} {\bf 23} (1999), 321-332.


\bibitem{Kamada}
H. Kamada, Self-dual K\"ahler metrics of neutral signature on complex surfaces,
PhD thesis, Tohoku University (2002).

\bibitem{kamada2} H. Kamada, Compact scalar-flat indefinite K\"ahler surfaces
with Hamiltonian $S^1$-symmetry, {\it Comm. Math. Phys.} {\bf 254} (2005),
23-44.

\bibitem{Kirch2006} K.-D.~Kirchberg, Integrability conditions for almost
Hermitian and almost K\"ahler $4$-manifolds, preprint, arXiv: math.DG/0605611.

\bibitem{kuiper} N. Kuiper, On conformally-flat spaces in the large, {\it Ann.
Math.} (2) {\bf 50} (1949), 916-924.


\bibitem{LawMat2} P. Law, Y. Matsushita, A spinor approach to Walker geometry,
{\it Comm. Math. Phys.} to appear, arXiv: math/0612804.


\bibitem{LB-M} C. LeBrun, L. Mason, Nonlinear gravitons, null geodesics and
holomorphic discs, {\it Duke Math. J.} {\bf 136} (2007), 205--273.



\bibitem{Mag1982}
M.~A.~Magid, {Indefinite Einstein hypersurfaces with nilpotent shape
operators}, \textit{Hokkaido Math. J.} \textbf{13} (1984), 241--250.


\bibitem{Matsushita1991}
Y.~Matsushita, {Fields of 2-planes and two kinds of almost complex structures
on compact 4-dimensional manifolds}, \textit{Math. Z.} {\bf 207} (1991),
281--291.

\bibitem{Mat}
Y.~Matsushita, { Walker $4$-manifolds with proper almost complex structure},
{\it J. Geom. Phys.} {\bf 55} (2005), 385-398.

\bibitem{LawMatsushita} Y. Matsushita, P. Law, Hitchin-Thorpe type inequalities
for pseudo-Riemannian 4-manifolds of metric signature (+ + - -), {\it Geom. Dedicata}  {\bf 87} (2001), 65-89.


\bibitem{Nak} I. Nakamura, Towards classification of non-K\"ahler complex
surfaces, {\it Sugaku Expos.}  {\bf 2} (1989), 209--229.


\bibitem{Nurowski-Przanowski1999}
P.~Nurowski, M.~Przanowski, {A four-dimensional example of a Ricci flat metric
admitting almost-K\"ahler non-K\"ahler structure}, {\it Classical Quantum
Gravity} {\bf 16}, no.\,3, (1999), L9-L13.

\bibitem{OV}
H. Ooguri, C.Vafa, { Self-duality and $N=2$ string magic}, {\it Modern Phys.
Lett. A} {\bf 5}, no.18 (1990), 1389-1398.



\bibitem{Ogr-Seki1998}
T.~Oguro, K.~Sekigawa, {Four-dimensional almost K\"ahler Einstein and
$\ast$-Einstein manifolds}, \textit{Geom. Dedicata} {\bf 69} (1998), 91--112.

\bibitem{Ols1978}
Z. Olszak, A note  on almost-K\"ahler manifolds, {\it Bull. Acad. Polon. Sci.
S\'er. Sci. Math. Astronom. Phys.}   {\bf 26} (1978), 139--141.


\bibitem{Ovando2004}
G.~Ovando, {Invariant pseudo K\"ahler metrics in dimension four}, {\it J. Lie
Theory} {\bf 16} (2006), 371--391.


\bibitem{P}
J. Petean, Indefinite K\"ahler-Einstein metrics on compact complex surfaces,
{\it Commun. Math. Phys.} {\bf 189} (1997), 227-235.




\bibitem{Sek}
K. Sekigawa, { On some compact Einstein almost K\"ahler manifolds}, {\it J.
Math. Soc. Japan} \textbf{36} (1987) 677--684.


\bibitem{tod} K. Tod, Indefinite conformally-ASD metrics on $S^2\times S^2$. In
{\it Further advances in twistor theory} (2001), vol. III, Chapman \& Hall/CRC, pp. 61-63.

\bibitem{Walker1950a}
A.~G.~Walker, {Canonical form for a Riemannian space with a parallel field of
null planes}, \textit{Quart. J. Math. Oxford (2)} {\bf 1}, 69--79 (1950).




























\end{thebibliography}
\end{document}